\documentclass[review]{elsarticle}
\usepackage{lineno,hyperref}
\usepackage{graphicx}
\modulolinenumbers[10]
\journal{OR Letters}

\usepackage{amssymb,amsmath,amsfonts}
\usepackage{algorithmic}
\usepackage{algorithm}

\newtheorem{theorem}{Theorem}%

\newtheorem{corollary}{Corollary}%

\begin{document}

\begin{frontmatter}

\title{A $4/3\cdot OPT+2/3$ approximation for big two-bar charts packing problem }

\author{Adil Erzin, Alexander Kononov}
\address{Sobolev Institute of Mathematics, Novosibirsk, Russia}
\author{Georgii Melidi}
\address{Sorbonne University, CNRS, Laboratoire d'informatique de Paris 6, France}
\author{Stepan Nazarenko}
\address{Novosibirsk State University, Novosibirsk, Russia}

\begin{abstract}
Two-Bar Charts Packing Problem is to pack $n$ two-bar charts (2-BCs) in a minimal number of unit-capacity bins. This problem generalizes
the strongly NP-hard Bin Packing Problem. We prove that the problem remains strongly NP-hard even if each 2-BC has at least one bar
higher than 1/2.  Next we consider the case when the first (or second) bar of each 2-BC is higher than 1/2 and show that the $O(n^2)$-time greedy algorithm with preliminary lexicographic ordering of 2-BCs constructs a packing of length at most $OPT+1$, where $OPT$ is optimum. Eventually, this result allowed us to present an $O(n^{2.5})$-time algorithm that constructs a packing of length at most $4/3\cdot OPT+2/3$ for the NP-hard case when each 2-BC has at least one bar higher than 1/2.
\end{abstract}


\end{frontmatter}

\linenumbers

\section{Introduction}

We consider the following resource-constrained scheduling problem. Given a set of jobs with two no-wait unit-execution-time operations. Both operations of the same job must be performed without delay. All jobs can run in parallel but each operation requires for its processing certain fixed amount of a common renewable resource. This scheduling problem was recently introduced in
\cite{Erzin20_1, Erzin21_2, Erzin21_3, Erzin21_4} and was formulated as a packing problem. Given a set of bar charts (2-BCs) consisting of two bars  and a fixed sequence of bins of unit capasity. Each bar must be packed into a bin so that the bars of each 2-BC do not change their order and they occupy adjacent bins. The goal is to minimize the number of used bins. We will use the term packing length to mean the number of used bins. The problem is called the Two-Bar Charts Packing Problem (2-BCPP) for the first time in \cite{Erzin21_2} and in subsequent papers. It is a generalization of the Bin Packing Problem (BPP) \cite{Johnson73} and a relaxation of the Two-Dimensional Vector Packing Problem (2-DVPP) \cite{Garey76,Kellerer03}.

In the Bin Packing Problem (BPP), items with given sizes must be packed in a minimal number of unit-capacity bins.
BPP is a particular case of 2-BCPP when each 2-BC consists of equal bars. The Bin Packing Problem is strongly NP-hard
and even the existence of a $(3/2-\varepsilon)$-approximation algorithm for the problem implies P=NP.
In the 70s and 80s of the last century, greedy algorithms were proposed, such as First-Fit, Best-Fit, First-Fit Decreasing etc \cite{Johnson73,Baker85,Yue91,Li97,Dosa07,Johnson85}. The main advantage of such algorithms is their low complexity.
Garey and Johnson \cite{Johnson85} proposed the Modified First-Fit Decreasing algorithm (MFFD) that has the best guaranteed accuracy among the known greedy algorithms. As shown in \cite{Yue95}, the packing length obtained by MFFD is bounded by $71/60\cdot OPT+1$, where $OPT$ is optimum. In their celebrated work, de la Vega and Lueker \cite{Vega1981} gave the first APTAS for  the Bin Packing Problem.
They showed that for any fixed $\varepsilon > 0,$ there exists a polynomial-time algorithm with asymptotic worst-case ratio not exceeding $1+\varepsilon.$ The result was consistently improved first by Karmarkar and Karp \cite{Karmarkar82} and then by Hoberg and Rothvoss \cite{Hoberg2017}. In the last paper mentioned above, the approximation ratio was improved to $OPT + O(\log OPT).$
For a detailed survey we refer the interested reader to \cite{Coffman2013, Christensen17}.

The Two-Dimensional Vector Packing Problem considers two attributes for each item and bin. The problem is to pack all items into
the minimal number of bins, considering both attributes of the bin's capacity limits.
The Two-Bar Charts Packing Problem is a relaxation of Two-Dimensional Vector Packing Problem. Indeed, every feasible solution
of 2-DVPP is a feasible solution of 2-BCPP.  At the same time, the value of the objective functions differs exactly by a factor of two.
In turn, any solution of 2-BCPP can be transformed into a solution of 2-DVPP, so that the number of occupied bins will increase at most by a factor of two. Thus, any $\rho$-approximation algorithm for 2-DVPP is a $2\rho$-approximation algorithm for 2-BCPP.
From approximation point of view, the Two-Dimensional Vector Packing Problem is harder than  the Bin Packing Problem.
Woeginger \cite{Woeginger1997} proved that there is no asymptotic polynomial time approximation scheme unless $P = NP$.
 The best known algorithm yields a $(3/2+\varepsilon)$-approximate solution for any $\varepsilon >0$ \cite{Bansal16}.
 See \cite{Christensen17} for a detailed review of approximations algorithms for 2-DVPP.

The Bar Charts Packing problem was first formulated in \cite{Erzin20_1} to optimize investments in the development of oil and gas fields.
Then the Two-Bar Charts Packing Problem was considered in \cite{Erzin21_2,Erzin21_3,Erzin21_4, Erzin22}.
Erzin et al. \cite{Erzin21_2} proposed a linear time algorithm for the 2-BCPP that gives a packing of length of at most $2OPT+1$,
where $OPT$ is the minimal packing length. A bar is called \emph{big} if it is higher than 1/2. A 2-BC is called \emph{big} if at least one of its bars is big. A particular case of the 2-BCPP where all 2-BCs are big was considered in  \cite{Erzin21_3,Erzin22}.
We denote this problem by 2-BCPP$^>$. An algorithm for finding a 3/2-approximation of optimal packing length was given
 in \cite{Erzin21_3}. This result was improved in \cite{Erzin22} where a 16/11-approximation algorithm for 2-BCPP$^>$ was presented.
 In the same paper, the authors proposed a 5/4-approximation algorithm if all the first (second) bars are big.

We note that despite the existence of approximate algorithms for various variants of 2-BCPP with big bars,
the question of their computational complexity remains open. In the next section we will prove that 2-BCPP$^>$ is
strongly NP-hard. In Section 3 we consider a greedy heuristic presented in  \cite{Erzin21_3,Erzin21_4}
and show that it finds a solution with a packing length at most $OPT+1$ for the case when all the first (second) bars are big.
Using this result we present $2/3(2OPT+1)$-approximation algorithm for 2-BCPP$^>$ in Section 4.


We end this section with exact definition of 2-BCPP$^>$.
Given a set $\cal S$ of $n$ big 2-BCs and $2n$ bins numbered from 1 to $2n.$
For each two-bar chart $i\in {\cal S}$ we denote the height of the first (left) and second (right) bar as $a_i$ and $b_i,$  respectively.
We assume that $\max\{a_i, b_i\} > \frac12$ for all $i\in {\cal S}.$ Let us denote by $[k, r]$ the set of integers $\{k, k+1, \ldots, r\}$.


\textbf{Definition 1.}
\emph{Packing} is a function $p: {\cal S} \rightarrow [1, 2n-1],$ which associates with each 2-BC $i \in {\cal S}$
an integer $p(i)$ corresponding to the bin number into which its first bar falls.
Let $B_k=\{i \in {\cal S}|p(i)=k\}.$ To avoid uncertainty we set $B_0= \emptyset.$
The packing is \emph{feasible} if the sum of the bar's heights that fall into each bin does not exceed 1,
i.e., $\sum_{i \in B_k} a_{i}+\sum_{i \in B_{k-1}} b_i \leq 1$ for all $k\in [1,2n]$.

As a result of a packing $p$, the first bar of 2-BC $i$ falls into the bin $p(i)$ and the second bar falls into the bin $p(i)+1$.
We will consider only feasible packings; therefore, the word ``feasible'' will be omitted further.

\textbf{Definition 2.}
The packing \emph{length} $L(p)$ is the number of bins in which at least one bar falls.


\textbf{The problem 2-BCPP$^>$ is to pack the big 2-BCs into the minimal number of bins.}


\section{NP-hardness of 2-BCPP$^>$}

In this section, we show that 2-BCPP$^>$ is strongly NP-hard. Our proof is based on a reduction from
Numerical 3-Dimensional Matching (SP16 in \cite{Garey79}).\\
{\large Numerical 3-Dimensional Matching} \\
 {\it Instance:} Three disjoint sets of positive integers $X=\{x_1,\ldots, x_r\},$ 
$Y=\{y_1,\ldots, y_r\},$ and $Z=\{z_1,\ldots, z_r\}$, consisting of $r$ elements each,
and the integer $b > 0$. \\
{\it Question}: Can $X\bigcup Y\bigcup Z$ be partioned into $r$ disjoint sets $A_1, A_2, \ldots, A_r$ such that each $A_i,$
$i=1,\ldots,r,$ contains one element from each of $X,$ $Y,$ and $Z,$ and $\sum_{e_j \in A_i} e_j = b.$


\begin{theorem}
The 2-BCPP$^>$ is strongly NP-hard.
\end{theorem}

\begin{figure}
\centering
\includegraphics[bb= 0 0 700 250, clip, scale=0.5]{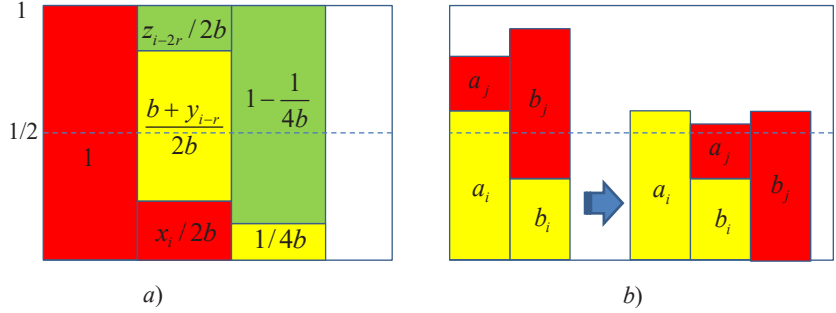}
\caption{\emph{a}) Illustration for the proof of Theorem 1; \emph{b}) Illustration for the proof of Theorem 3.} \label{fig1}
\end{figure}

\textbf{Proof.}
Let $\{x_1,\ldots, x_r\},$ $\{y_1,\ldots, y_r\},$  $\{z_1,\ldots, z_r\}$ and the integer $b > 0$ be an arbitrary instance $I$ of the Numerical 3-Dimensional Matching problem. Without loss of generality we assume that $\sum_{i=1}^{r}(x_i+y_i+z_i) = rb.$ We construct an
instance $I'$ of 2-BCPP$^>$ with $3r$ 2-BCs as follows:
\begin{itemize}
\item Red 2-BCs $i,$ $1\leq i \leq r,$ $a_i = 1,$ $b_i=\frac{x_i}{2b}$,
\item Yellow 2-BCs $i,$ $r+1\leq i \leq 2r,$ $a_i = \frac{b+y_{(i-r)}}{2b},$ $b_i=\frac{1}{4b}$,
\item Green 2-BCs $i,$ $2r+1\leq i \leq 3r,$ $a_i = \frac{z_{(i-2r)}}{2b},$ $b_i=1-\frac{1}{4b}.$
\end{itemize}
We note that $a_i > \frac12$ for $1\leq i \leq 2r$ and $b_i > \frac12$ for  $2r+1\leq i \leq 3r.$
Thus, we definitely get an instance of 2-BCPP$^>$.
We claim that the answer to instance $I$ is affirmative if and only
if there is a feasible solution of instance $I'$ which uses at most $3r$ bins.

First, let the required partition exist and $x_i + y_j + z_k = b$ for some $i, j, k.$
Then the three 2-BCs $i,$ $r+j,$ and $2r+k$ can be packed into three bins as shown in Figure~\ref{fig1}$a$.
Hence, there is a feasible solution of instance $I'$ which uses at most $3r$ bins.

Next, let $p$ be a packing such that $L(p) \leq 3r.$ Since all 2-BCs are big each bin contains at most three bars.
Indeed, let some bin contain more than three bars. No more than one big bar can be placed in each bin.
Therefore, all other bars in the bin should be small. But then either the left or right bin contains two big bars.
We got a contradiction with the feasibility of packing.

Since $a_i = 1$ for every $1\leq i \leq r$, $r$ bins should be occupied by exactly one bar, namely, the left bar of the red 2-BC.
Therefore, of the remaining $2r$ bins, $r$ bins contain exactly 3 bars in each bin, and $r$ bins contain exactly 2 bars in each bin.
Since  $b_i=1-\frac{1}{4b}$ for every $2r+1\leq i \leq 3r,$ and $\min_{i=1,\ldots,r}\{x_i, y_i, z_i\}\geq 1,$
the right bar of green 2-BC can be placed in the same bin only with the right bar of yellow 2-BC.
Such pairs will take up another $r$ of bins. Hence, the right bars of the red 2-BC,
the left bars of the yellow 2-BC and the left bars of the green 2-BC should be distributed in the remaining $r$ bins of three each.
We have $$\sum_{i=1}^{r}b_i +\sum_{i=r+1}^{2r}a_i + \sum_{i=2r+1}^{3r}a_i = \sum_{i=1}^r \frac{b+x_i+y_i+z_i}{2b} = r.$$
Therefore, the remaining $3r$ bars with total size of $r$ are located in $r$ unit-capacity bins.
It is possible if and only if the instance $I$ of  Numerical 3-Dimensional Matching has the required partition. \qed

\section{Linearly ordered packing}
In this section, we consider the particular case of the 2-BCPP,
in which we will impose an additional restriction on the set of feasible solution.
We say that a packing $p$ is {\it linearly ordered} if for any two 2-BCs $x, y \in \cal S$ we have $p(x) \ne p(y).$
We note that with a linearly ordered packing, each bin contains no more than two bars.
In this section, we present a greedy algorithm that finds an almost optimal  linearly ordered packing.
At the end of the section, as a consequence, we will get an approximation with an absolute error bound for the
2-BCPP$^>$ problem when the first bar of each  2-BC is big.

Let $p: {\cal S} \rightarrow [1, 2n-1],$ be a linearly ordered packing and $x, y \in {\cal S}$ such that
$p(x)=i$ and $p(y)=i+1$ for some $i \in [1,2n-2].$ We refer to $x$ as the left neighbour of $y$
and to $y$ as the right neighbour of $x$, respectively, and we write $r(x)=y$, i.e., $r(x)$ denotes the right neighbour of  $x$.
We call the maximal sequence of 2-BCs, in which each 2-BC except the first one has a left neighbour, a \textit{chain}.
Let $\lambda(p)$ be a number of chains in the packing $p$, then $L(p)=n + \lambda(p).$

We adopt the greedy algorithm $GALO$ presented in \cite{Erzin21_2,Erzin21_4} for this particular case.
Algorithm $GALO$ works as follows.
Sort all 2-BCs lexicographically in non-increasing order of bar's height. Starting from the first bin,
repeat the following procedure until all 2-BCs are packed. Place first 2-BC from the list (if any) that can be packed
into the current bin. Remove this 2-BC from the list and move to the next bin.

Algorithm $GALO$ performs at most $2n-1$ iterations. At each iteration, the search for a suitable 2-BC
requires no more than $O(n)$ time. Thus, the running time of algorithm $GALO$ is $O(n^2)$.

Denote the packing constructed by algorithm $GALO$ as $p^g.$

\begin{theorem}
Let $p^*$ be an optimal linearly ordered packing,  then $L(p^g) \leq L(p^*)+1$.
\end{theorem}

\textbf{Proof.} After the first step of algorithm $GALO$, all 2-BCs are ordered lexicographically in a non-increasing manner.
For two 2-BCs: $x$ and $y$, we will write  $x\leftarrow y$ if $b_x + a_y \leq 1.$

Let
$
i(x)=\min\{y\in {\cal S}\mid x\leftarrow y\},\ x\in {\cal S}.
$
If for some $x\in {\cal S}$ there is no any 2-BC $y\in {\cal S}$ such that $b_x + a_y \leq 1$, then set (for definiteness) $i(x)=n+1$.
Note that for any $x, y \in {\cal S}$ such that
$y \geq i(x)$ we have $x\leftarrow y$.

We call $x$ in the packing an \textit{unfortunate} 2-BC if it does not have a right neighbour.
Let us number the chains in the packing $p^g$ in the order of their appearance in $p^g$, and refer to the $i$-th chain as $C_i$.
It is clear that the last 2-BC in each chain is unfortunate  2-BC.
We will prove that all but one of the unfortunate 2-BCs in the packing $p^g$
can be put into one-to-one correspondence with the unfortunate 2-BCs in the optimal packing $p^*$.

Let $x_1 \in C_1$ be the first unfortunate 2-BC in the packing $p^g$.

\begin{enumerate}
  \item If $i(x_1)=n+1$,  then none of 2-BC can be the right neighbour of $x_1$.
   Therefore, $x_1$ also does not have a right neighbour in the optimal packing $p^*$.
  \item Let $i(x_1)\leq n$.  Then all 2-BCs in the packing $p^g$ from $i(x_1)$ up to $n$ belong to $C_1$.
  Denote the set of their left neighbours by $U_1$.  Let
$$
i(U_1)=\min\{i(y) \mid y\in U_1\}.
$$
Note that $i(U_1)\leq i(x_1)$.

\begin{itemize}
\item[2.0] If $i(x_1)=1$,  then $U_1 \equiv  {\cal S}$ and,  consequently,  all 2-BCs except $x_1$ have a right neighbour.
Since there should be at least one unfortunate 2-BC (one chain) in the optimal packing, the resulting solution is optimal.
\end{itemize}
Let $i(x_1)>1$.
\begin{itemize}
\item[2.1] If $i(U_1)=i(x_1)$, then there is at least one 2-BC in $U_1\cup\{x_1\}$ which in $p^*$  does not have a 1-union on the right.  Indeed, 
we have $ \left| U_1 \cup\{x_1\}\right| > \left|\{i(x_1),\dots,n\}\right|$ and at least one 2-BC from the set  $U_1 \cup\{x_1\}$  will not get a right neighbour. Choose such an arbitrary unfortunate 2-BC $x_1'\in U_1\cup \{x_1\}$ in $p^*$ and put it in accordance with $x_1$.

\item[2.2]  Let $i(U_1)<i(x_1)$.  Choose such $z\in {\cal S}$ that $i(z)=i(U_1)$. Since $i(z) < i(x_1)$ and $r(z) \geq i(x_1)$, all 2-BCs from $i(z)$ to $i(x_1)-1$ are packed before $z$ and belong to $C_1.$  We set
$$
U_1=U_1\cup\{u\in S \mid u\leftarrow y,\ i(z)\leq y\leq i(x_1)-1\}.
$$
Repeating the above reasoning for $z$, we again find ourselves in one of the three cases described above.
In case 2.0, the solution is optimal.  In case 2.1, in the optimal packing $p^*$, there is an unfortunate 2-BC in the set $U_1\cup\{z\}$,  which corresponds to $x_1$. If we are in case 2.2, we will find a new element $z'$ such that $i(z') < i(z)$ and $r(z') \geq i(z)$ and repeat
our reasoning. Each time in case 2.2,  the value of $i(U_1)$ decreases, and, eventually, at some iteration, we get into case 2.0 or 2.1.
\end{itemize}

Suppose we have matched $k-1$ unfortunate 2-BCs in the packing $p^g$ with unfortunate 2-BCs in $p^*$. By construction, all these 2-BCs belong to the set $U_{k-1}\cup\{x_1,x_2,\ldots, x_{k-1}\}$.

\item Let $x_k \in C_k$ be the $k$th unfortunate 2-BC (which has no right neighbour in the packing $p^g$).  Let
$$
r_k=\min\{i\mid i\notin C_1\cup\dots\cup C_{k-1}\}.
$$
Algorithm $GALO$ starts a new chain with $r_k$.
Note that all 2-BCs with numbers greater than or equal to $i(U_{k-1})$  and less than $r_k$  have already been used in the previous chains.

\begin{itemize}
\item[3.0] If $i(x_k)\leq r_k$,  then all 2-BCs from $r_k+1$ to $n$ either belong to the first $k-1$ chains,  or have a left neighbour in the chain $C_k$. But,  as noted above,  all 2-BCs with numbers less than $r_k$ are also used in the previous chains.  Therefore,  the chain $C_k$  is the last one in the packing $p^g$.  Hence,  in $p^g$,  for all unfortunate 2-BCs except $x_k$,  a one-to-one correspondence with unfortunate 2-BCs in the packing $p^*$ is found,  so $L(p^g)\leq L(p^*)+1$.
\end{itemize}
Let $i(x_k)>r_k$.

Set $U_k$ as follows. If $i(U_{k-1})\leq i(x_k)$,  then $U_k=U_{k-1}$. Let $i(U_{k-1})>i(x_k)$. In this case,  all 2-BCs in the packing $p^g$ from $i(x_k)$ to $n$ have the left neighbours.  Denote this set by $U_k$.  Let
$$
i(U_k)=\min\{i(y)\mid y\in U_k\}.
$$
Note that $i(U_k)\leq i(x_k)$.

\begin{itemize}
\item[3.1] If $i(U_k)=i(x_k)$, then there is at least one 2-BC in $U_k\cup\{x_k\}$ which in $p^*$ does not have a right neighbour,  not matching the previously selected $\{x_1',\dots, x_{k-1}'\}$.  Similarly to the case 2.1,  we select such an arbitrary 2-BC $x_k'$ and put it in correspondence with $x_k$.
 \item[3.2] Let $i(U_k)<i(x_k)$.  Choose such 2-BC $z$ that $i(z)=i(U_k)$. Since $i(z)<i(x_k)$,  all 2-BCs from $i(z)$ to $i(x_k)-1$ are already packed before $z$.  We set
$$
U_k=U_k\cup\{u\in S \mid u\leftarrow y,\ i(z)\leq y\leq i(x_k)-1\}.
$$
Repeating the reasoning for the 2-BC $z$, we again fall into one of the three cases described above.
Each time when we go back to case 3.2,  the value of $i(U_k)$ decreases, and, eventually, at some iteration, we get into case either 3.0 or 3.1.


Thus, at the $k$th step, we either find a new correspondence between the unfortunate 2-BCs in the packings $p^g$ and $p^*$, or observe that the chain under consideration is the last one in $p^g$ and there is one-to-one correspondence with unfortunate 2-BCs in the optimal packing $p^*$ for all unfortunate 2-BCs in $p^g$, except $x_k$.  Hence, we finally get the statement of the theorem.  
\qed
\end{itemize}
\end{enumerate}


 Let  $\cal S$ be a set of $n$ 2-BCs such that $a_j > \frac12$ for all $j \in [1,n].$
 Since for any 2-BCs $x, y \in {\cal S}$ we have $a_x + a_y > 1$, any feasible packing of  $\cal S$ is linearly ordered.
 So, Theorem 2 implies the following result.

\begin{corollary}
If the first bar of each 2-BC is big, then $L(p^g) \leq OPT+1,$ where $OPT$ is the length of optimal packing.
\end{corollary}

\section{A $4/3\cdot OPT+2/3$ approximation for 2-BCPP$^>$}
In this section, we consider the case when all 2-BCs are big (at least one bar of each 2-BC is big).
Let $p: {\cal S} \rightarrow [1, 2n-1],$ be a packing of $ {\cal S}.$ We say that two 2-BCs $i, j \in \cal S$
form a pair if $p(i)=p(j).$ We note that there are no pairs of 2-BCs in a linearly ordered packing.

First, we present an algorithm that finds a packing in which there is the maximum possible number of pairs.
Algorithm Matching builds a graph $G=(V,E)$ in which the vertices are the images of the 2-BCs.
The edge $(i,j)\in E$ if both 2-BCs $i$ and $j$ can be packed in the two consecutive bins, i.e.,
if $a_i + a_j \leq 1$ and $b_i + b_j \leq 1.$ After that the algorithm finds a maximal matching $M$ in $G$  \cite{Xie18}.
This matching determines the packing $p^m$ as follows. Each pair of 2-BCs $i$ and $j$ such that $(i,j)\in M$
is packed  in the two consecutive bins. Each of the remaining 2-BCs is assigned to the bins in any way.
The running time of Algorithm Matching is evidently dominated by finding a maximal matching and
can be estimated by $O(n^{2.5})$ time \cite{Xie18}.

Now we are ready to present an approximation algorithm for 2-BCPP$^>$.

\begin{algorithm}[H]
\begin{algorithmic}[1]

\STATE Find a linearly ordered packing $p^g$ by algorithm $GALO$
\STATE Find a packing $p^m$ by Algorithm Matching
\STATE Take the best of two solutions, say $p^a.$

\end{algorithmic}
\caption{Algorithm App}\label{approx_algorithm}
\end{algorithm}

\begin{theorem} If all 2-BCs are big, then Algorithm $App$ finds a packing of length at most $4/3\cdot OPT+2/3$ in $O(n^{2.5})$
time, where $OPT$ is the length of optimal packing.  This performance guarantee is tight.
\end{theorem}

\textbf{Proof.} Let $p^*$ be an optimal packing, i.e. $OPT = L(p^*).$ Since any two big bars cannot occupy the same bin,
we have $OPT \geq n.$ Denote by $k_2$ the number of pairs of 2-BCs in $p^*$.
Each pair consists of two 2-BCs and occupies two bins. Then the number of remaining 2-BCs is equal to $n-2k_2 \doteq k_1.$

Let $\mu$ be the number of pairs obtained by Algorithm Matching. Then, $L(p^m) \leq 2\mu + 2k_1.$
Since Algorithm Matching finds a packing with the maximum possible number of pairs, we have $\mu \geq  k_2.$
Thus, $L(p^m) \leq 2n - 2\mu = 2k_1+4k_2-2\mu \leq 2k_1 + 2k_2 = k_1 + n \leq OPT + k_1.$

Let $i, j \in \cal S$ be a pair in  $p^*$. Without lost of generality, we assume that $a_i > \frac12$ and $b_j > \frac12.$
Then, $a_j + b_i \leq 1.$ We will push the two bins $p(i)$ and $p(i)+1$ apart and place the left bar of $j$ and right bar of $i$ in a new bin between them.  We get a new feasible packing, the length of which will increase by 1 (Figure~\ref{fig1}$b$).
Having done this for all pairs in $p^*$ we get a linearly ordered packing with a packing length $OPT + k_2.$
Then Theorem 2 implies $L(p^g) \leq OPT + k_2 +1.$

Finally, we get $$L(p^a) \leq OPT + \min\{k_1,k_2+1\}  \leq OPT + \frac{k_1}{3} + \frac{2(k_2+1)}{3}$$
$$= OPT + \frac{n}{3} + \frac23 \leq \frac43 OPT + \frac23.$$

\begin{figure}
\centering
\includegraphics[bb= 0 0 500 250, clip, scale=0.5]{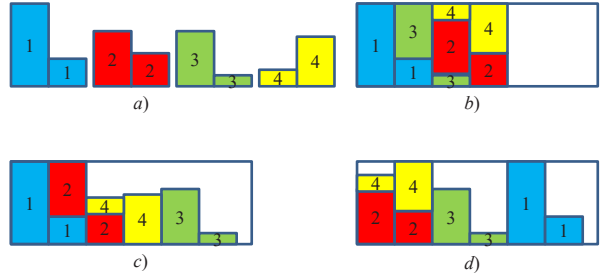}
\caption {\emph{a}) Lexicographically ordered set of 2-BCs; \emph{b}) Optimal packing; \emph{c}) Packing $p^g$ built by $GALO$; \emph{d}) Packing $p^m$ built by Matching.}
\label{fig4}
\end{figure}
Fig. \ref{fig4} shows that the performance guarantee obtained is tight.
In this example, the optimal packing length is 4, while the length of packings $p^g$ and $p^m$ is 6.
Then $L(p^a) =\frac{4}{3}OPT+\frac{2}{3}=6$.

Time complexity of the algorithm $GALO$ is $O(n^2)$. Maximal matching can be found with $O(n^{2.5})$ time complexity \cite{Xie18}. \qed

\section{Conclusion}
For the problem 2-BCPP$^>$ of packing big two-bar charts (the height of at least one bar is more than 1/2)
in a minimal number of unit-capacity bins:\\
(i) we prove that it is strongly NP-hard;\\
(ii) for the case when the first (second) bar of each 2-BC is big, we prove that algorithm $GALO$
yields a packing of length at most $OPT+1$ in $O(n^2)$ time;\\
(iii) we propose an algorithm which finds a packing for 2-BCPP$^>$ of length no more than $4/3\cdot OPT+2/3$ in $O(n^{2.5})$ time.

\paragraph{Acknowledgement} The study was carried out within the framework of the state contract of the Sobolev Institute of Mathematics (project FWNF-2022-0019).


\bibliography{ErKonMelNaz}

\end{document}